\documentclass[12pt]{article}
\usepackage{amsmath,amscd,amssymb,amsfonts,amsthm,diagrams}

\title{Dual Forms on Supermanifolds and Cartan Calculus}

\newcommand{\myaddress}{\begin{center}
\small Department of Mathematics\\ \small University of
Manchester Institute of Science and Technology (UMIST)\\ \small
PO Box 88, Manchester M60 1QD, England\\ \small {\tt
theodore.voronov@umist.ac.uk}
\end{center}}

\author{ Theodore Voronov}

\date{\myaddress }

\renewcommand{\leq}{\leqslant}
\renewcommand{\geq}{\geqslant}

\newcommand{\cowedge}{{\lrcorner}}
\DeclareMathOperator{\Ber}{Ber} \DeclareMathOperator{\Vol}{Vol}
\DeclareMathOperator{\E}{\Lambda}
\DeclareMathOperator{\EE}{{\boldsymbol\Lambda}}
\DeclareMathOperator{\OO}{{\boldsymbol \Omega}}
\DeclareMathOperator{\GL}{GL}
\DeclareMathOperator{\Cliff}{Cliff}
\DeclareMathOperator{\Mat}{Mat}

\DeclareMathOperator{\Vect}{Vect}
\newcommand{\der}[2]{{\frac{\partial {#1}}{\partial {#2}}}}
\newcommand{\lder}[2]{{\partial {#1}/\partial {#2}}}
\newcommand{\dder}[3]{{\frac{\partial^2 {#1}}{\partial {#2}\partial {#3}}}}
\newcommand{\da}{{\partial_A}}
\newcommand{\R}[1]{{\mathbb R}^{#1}}
\renewcommand{\L}{{\cal L}}
\newcommand{\Z}{{\mathbb Z_{2}}}
\newcommand{\dlim}{{\mathop{\varinjlim}\limits_{N,M}}}
\newcommand{\D}{{\,\,}{\bar{\smash{\!\!\mathit d}}}}
\newcommand{\DD}{{}{\bar{\smash{\delta}}}}
\newcommand{\bdot}{{\boldsymbol \cdot}}

\renewcommand{\a}{\alpha}
\renewcommand{\b}{\beta}

\renewcommand{\d}{\delta}
\newcommand{\e}{\varepsilon}

\newcommand{\x}{\xi}
\renewcommand{\o}{\omega}
\renewcommand{\O}{{\Omega}}
\newcommand{\s}{\sigma}
\renewcommand{\t}{\tau}
\renewcommand{\P}{\Pi}

\newtheorem{thm}{Theorem}[section]

\newtheorem{co}{Corollary}[section]

\theoremstyle{definition}
\newtheorem{de}{Definition}[section]
\newtheorem{ex}{Example}[section]

\newtheorem{rem}{Remark}[section]

\newcommand{\vfa}{{v_F{}^A}}
\newcommand{\vga}{{v_G{}^A}}
\newcommand{\vfb}{{v_F{}^B}}
\newcommand{\vgb}{{v_G{}^B}}

\newcommand{\pak}{{p_A{}^K}}
\newcommand{\pal}{{p_A{}^L}}

\newcommand{\pbk}{{p_B{}^K}}
\newcommand{\pdk}{{p_D{}^K}}
\newcommand{\pbl}{{p_B{}^L}}
\newcommand{\pdl}{{p_D{}^L}}
\newcommand{\pcl}{{p_C{}^L}}
\newcommand{\pck}{{p_C{}^K}}

\newcommand{\pakstar}{{p_A{}^{K^*}}}
\newcommand{\pbkstar}{{p_B{}^{K^*}}}
\newcommand{\wkstar}{{w^{K^*}}}

\newcommand{\wfk}{{w_F{}^K}}
\newcommand{\wfl}{{w_F{}^L}}

\newcommand{\wgk}{{w_G{}^K}}
\newcommand{\wgl}{{w_G{}^L}}

\newcommand{\vra}{{v_r{}^A}}
\newcommand{\wra}{{w_r{}^A}}
\newcommand{\wrk}{{w_r{}^K}}

\newcommand{\wk}{{w^K}}
\newcommand{\ua}{{u^A}}
\newcommand{\xa}{{x^A}}
\newcommand{\xb}{{x^B}}
\newcommand{\Xa}{{X^A}}
\newcommand{\Xb}{{X^B}}

\newcommand{\dab}{{\d_A{}^B}}

\newcommand{\kt}{{\tilde K}}
\newcommand{\lt}{{\tilde L}}
\newcommand{\at}{{\tilde A}}
\newcommand{\bt}{{\tilde B}}
\newcommand{\ct}{{\tilde C}}
\newcommand{\dt}{{\tilde D}}
\newcommand{\ft}{{\tilde F}}
\newcommand{\gt}{{\tilde G}}
\newcommand{\ut}{{\tilde u}}
\newcommand{\vt}{{\tilde v}}
\newcommand{\xt}{{\tilde X}}
\newcommand{\aat}{{\tilde \alpha}}
\newcommand{\bbt}{{\tilde \beta}}

\newcommand{\pone}{{p_1}}
\newcommand{\ptwo}{{p_2}}
\newcommand{\woo}{{w_{11}}}
\newcommand{\wot}{{w_{12}}}
\newcommand{\wto}{{w_{21}}}
\newcommand{\wtt}{{w_{22}}}

\newcommand{\ppone}[1]{{p_{#1}^{p+1}}}
\newcommand{\pptwo}[1]{{p_{#1}^{p+2}}}
\newcommand{\paone}{{\ppone{A}}}
\newcommand{\pbone}{{\ppone{A}}}
\newcommand{\pcone}{{\ppone{C}}}
\newcommand{\pdone}{{\ppone{D}}}

\newcommand{\patwo}{{\pptwo{A}}}
\newcommand{\pbtwo}{{\pptwo{B}}}
\newcommand{\pctwo}{{\pptwo{C}}}
\newcommand{\pdtwo}{{\pptwo{D}}}

\newcommand{\wwone}[1]{{w^{#1}_{r+1}}}
\newcommand{\wwtwo}[1]{{w^{#1}_{r+2}}}
\newcommand{\wkone}{{\wwone{K}}}
\newcommand{\wktwo}{{\wwtwo{K}}}
\newcommand{\wlone}{{\wwone{L}}}
\newcommand{\wltwo}{{\wwtwo{L}}}
\newcommand{\wtwo}{{w^{p+2}}}
\newcommand{\wone}{{w^{p+1}}}

\newcommand{\pnone}{{p^{n+1}}}
\newcommand{\wnone}{{w^{n+1}}}
\newcommand{\panone}{{p_A{}^{n+1}}}
\newcommand{\pbnone}{{p_B{}^{n+1}}}
\newcommand{\wfnone}{{w_F{}^{n+1}}}

\newcommand{\vvone}[1]{{v^{#1}_{r+1}}}
\newcommand{\vaone}{{\vvone{A}}}
\newcommand{\vbone}{{\vvone{B}}}

\begin{document}

\maketitle

\begin{abstract}
\noindent
The complex of ``stable forms" on supermanifolds is studied.
Stable forms on $M$ are represented by
certain  Lagrangians of  ``copaths" (formal systems of equations, which may
or may not specify actual surfaces) on $M\times\mathbb R^D$.
Changes of $D$ give rise to
stability isomorphisms. The {Cartan-de~Rham} complex made of stable
forms extends both in positive and negative degree and its
positive half is isomorphic to the complex of forms defined as
Lagrangians of paths. Considering the negative half is necessary, in particular, for
homotopy invariance.

We introduce analogs of exterior multiplication by covectors and
of contraction with vectors. We find (anti)commutation relations for them.
An analog of Cartan's homotopy identity is proved. Before stabilization it contains
a stability operator $\s$.
\end{abstract}

\section*{Introduction}

The crucial difference of ``exterior algebra" in the super case
from the usual case is that the analog of the ``top exterior
power" for a $\Z$-graded vector space cannot be obtained by
tensor operations. This is because the determinant in the super
case (the Berezinian) is not a polynomial expression, but a
fraction whose numerator and denominator separately are not
multiplicative. Thus the space $\Ber V$ (which corresponds to the
usual $\det V$) enters independently of  the
``naive" generalization of exterior multiplication by the sign
rule. A complete theory of ``exterior forms" has to be built upon
the Berezinian from the beginning. This fact has far reaching
consequences.

``Naive" differential forms on a supermanifold $M^{n|m}$  are, of
course,
(locally)
polynomials in $dx^A$, where $x^A$ are coordinates.
Experts know that there are
two possible conventions for the parity and  commutation relations
for the differentials (see~\cite{manin}). According to one of
them, $dx^A$ is assigned the same parity as $x^A$ and the
differentials anticommute: the flip of $dx^A$ and $dx^B$ results in the factor
$-(-1)^{\tilde A \tilde B}$. The other convention assigns to
$dx^A$ the parity opposite to that of $x^A$ and the differentials
are regarded as commuting variables.
We shall refer to them in the sequel as to skew-commutative and
commutative conventions, respectively.

From the viewpoint of integration, the fatal drawback of such naive
forms is that they can't be integrated over $M=M^{n|m}$ (unless $m=0$).
Because of that, some remedies were suggested.

Bernstein and Leites~\cite{berl:int}
defined ``integral forms" as tensor products of multivector fields
with Berezin volume forms. This permitted integration over
$M^{n|m}$ and an analog of Gauss-Ostrogradsky formula. If we are
integration-minded, we expect that the correct forms on
supermanifolds should be graded by super dimensions $r|s$ (dimensions of
surfaces or chains over which a form can be
integrated). Thus, integral forms should correspond to
``$r|m$-forms" ($s=m$) and volume forms to ``$n|m$-forms".
Naive differential forms from this point of view correspond to
``$r|0$-forms" ($s=0$.) What about other values of $r$, $s$?

For  {\it non-polynomial} functions of $dx^A$ (with the
commutative convention) Bernstein and
Leites~\cite{berl:pdf} showed that they also can be integrated
over $M^{n|m}$ provided they sufficiently rapidly decrease in
$d\x^{\mu}$, where $\x^{\mu}$ are odd coordinates.
Such ``pseudodifferential forms" are very beautiful.
However, since they do not
have any grading (and, in fact, are good for integration
only for a particular type of orientation and not good for others,
see~\cite{tv:git}) they do not
solve the problem.

A crucial step towards the theory of ``$r|s$-forms"
was made  by
A.S.~Schwarz, M.A.~Baranov, A.V.~Gajduk,
O.M.~Khudaverdian and A.A.~Rosly in the beginning of 1980-s
and was motivated by quantum field theory.
They based
their investigation of the ``objects of integration" on
supermanifolds directly on the notion of  Berezinian and
studied  Lagrangians of
parameterized surfaces $\Gamma: I^{r|s}\to M^{n|m}$ which
induce volume forms on $r|s$-dimensional space $U^{r|s}$.
They are called  {\it densities}. The
key result was the concept of  ``closedness" of a
density~\cite{ass:bar,ghs,hov:viniti}: a density is said to be
{\it closed} if the corresponding action is identically
stationary. (On ordinary manifolds, for densities corresponding to
closed forms this follows from the Stokes' formula.)

As the author
discovered, this notion
of ``closedness" precisely follows from a certain construction of a
differential in terms of variational derivatives. Densities,
initially defined only for embedded
surfaces (hence $0\leq r\leq n$, $0\leq s\leq m$), should be
replaced by more general ``covariant Lagrangians", for which
$r\geq 0$ can exceed $n$, and a certain system of differential
equations with respect to the components of tangent vectors is
imposed upon Lagrangians. Roughly speaking, this system (see
Eq.~(\ref{eqs})
below) is a nontrivial analog of multilinearity/skew symmetry
property of the usual exterior forms. (The
odd-odd part of the system amazingly coincides with the equations introduced
by F.~John~\cite{john} and Gelfand-Shapiro-Gindikin-Graev
(see~\cite{gel:ggg}) for the description of the image of
Radon-like transforms  in integral geometry.)
The theory of $r|s$-forms in this sense was developed by the author with
A.V.~Zori{\'c}~\cite{tv:compl,tv:pdf,tv:bord, tv:coh} and the
author~\cite{tv:git}. The differential has degree $+1$, so
$r|s$-forms are mapped to $r+1|s$-forms. The complex obtained in
this way possesses all natural properties of the usual
Cartan-de~Rham complex like functoriality in a suitable category,
Stokes' formula and homotopy invariance, and also has some
similarity with extraordinary cohomology (an analog of the
Atiyah-Hirzebruch spectral sequence), see~\cite{tv:git}. For
$s=0$, it naturally incorporates the ``naive" generalization of
differential forms. For $s=m$ and $r\geq 0$ it also incorporated
integral forms of Bernstein and Leites.

However, an {\it ad hoc} augmentation of the complex
had to be introduced~\cite{tv:git} to achieve homotopy invariance.
The existence of Bernstein-Leites integral forms of negative degree
has also hinted to  ``hidden" $r|s$-forms with $r<0$.

Such objects were indeed discovered in~\cite{tv:dual}. Together with
forms considered in~\cite{tv:git} they give a
desired de Rham complex stretching  both in positive and negative directions.

The solution is based on the idea of a {\it dual form}~\cite{tv:dual}
(important results were
independently obtained in~\cite{hov:bv}).
Geometrically, dual forms  are Lagrangians of surfaces specified by
maps $M^{n|m}\supset U^{n|m}\to\mathbb R^{p|q}$
(copaths) rather than  maps $I^{r|s}\to M^{n|m}$ (paths). To define a
complex, dual forms are not enough. One has to introduce new
independent parameters and to allow to increase their number.
An intermediate product is labeled ``mixed form".
A whole
bunch of isomorphisms enters the stage, and the final picture
is the result of a stabilization (see~\cite{tv:dual} and subsection~\ref{algebra:st}
below).
(Geometrically, one  gets a sort of
``virtual surfaces", which can have both negative and positive
dimension.)

\bigskip
In the current paper we develop the algebraic and differential theory of
{\it stable forms} (the
unified complex). We do not touch integration.
The main result of the paper is an analog of Cartan calculus that
includes module structures for forms and the relation between the
differential, Lie derivative and a ``contraction operator" with a
vector field (which is defined in this paper). All results are
new. They will be used to study the homotopy properties of stable
forms and the de~Rham cohomology of supermanifolds.

The paper is organized as follows.

In Section~\ref{algebra} we define dual and mixed forms on a superspace $V$,
the stability isomorphisms and isomorphisms with forms
considered in~\cite{tv:git}. Operators $e(\a)$ and
$e(v)$ are introduced, where $u\in V$, $\a\in V^*$.
We prove that they are stable (commute with the stability isomorphisms)
and relate them with operators on forms of~\cite{tv:git}
(Theorem~\ref{stability}).
Then we find
the relations that they obey. We get a ``skew-commutative" version of a Clifford algebra
involving a stability operator $\s$ as an additional central element
(Theorem~\ref{commut}).
As a corollary, we obtain  module structures over the
exterior algebras $\E (V^*)$ and $\E (V)$ (the skew-commutative versions).

In Section~\ref{analiz} we consider the complex of stable forms
on a supermanifold $M$. We prove the Leibniz identity (=differential module structure)
for the multiplication by naive differential forms
$\o\in\O^{\bdot}(M)$
(Theorem~\ref{leibniz}).
Then we consider the Lie derivative for mixed forms.
We prove that the anticommutator of the differential and  the operator $e(X)$,
where $X$ is a
vector field, equals the Lie derivative multiplied by the operator $\s$
(Theorem~\ref{cart}).
It immediately implies a ``Cartan's homotopy identity" for stable
forms.

The results are discussed in Section~\ref{discuss}.

We mainly follow the notation and terminology of the
book~\cite{tv:git}.

\smallskip
{\bf Acknowledgements:} Questions related to the topic of this paper were
discussed at various times with J.N.~Bernstein, O.M.~Khudaverdian
and A.~Belopolsky. I am very much grateful to them.

\section{Algebraic theory}\label{algebra}
\subsection{Construction of forms. Stability isomorphisms}
\label{algebra:st}
Consider a superspace $V$ over $\mathbb R$ of dimension $\dim
V=n|m$.  We
identify vector superspaces with the corresponding supermanifolds.
By $\Vol V:=\Ber V^*$ we denote the space of volume forms on $V$.
In the following we consider functions whose arguments are vectors or covectors.
Components of vectors are written as rows, components of covectors
as columns.

Recall the following definition.

\begin{de}[{\normalfont see~\cite{tv:compl,tv:coh,tv:git}}]\label{straight}
A {\it form} on $V$ of degree $r|s$ is a smooth
map $L:\underbrace{V\times\dots V}_{r}\times
\underbrace{\P V\times\dots \times\P V}_{s} \to \mathbb R$ satisfying
the following
conditions~\eqref{bers} and~\eqref{eqs}:
\begin{align}
  &L(gv)= L(v) \Ber g,\label{bers}\\
\intertext{for all $g\in \GL(r|s)$ and}
  &\dder{L}{\vfa}{\vgb}
  +(-1)^{\ft\gt+(\ft+\gt)\bt}\dder{L}{\vga}{\vfb}=0\label{eqs}.
\end{align}
In our notation the argument of the function $L$ is
written as a matrix $v=(\vfa)$  whose rows $v_F$ are vectors (written
in components). The condition~\eqref{bers} implies that $L(v)$ is defined only
if
odd vectors $v_K$, $\kt=1$, are linearly independent. Hence $0\leq s \leq
m$, while $r\geq 0$ can be arbitrary.
\end{de}

Though this definition provides no efficient description of forms,
such a description can be given in special cases (corresponding to naive differential forms
and to Bernstein-Leites integral forms) and in other cases
various examples can be provided. See~\cite{tv:git}. In
particular, if $m>0$, for $s\neq m$ there are nonzero forms with
$r>n$. We shall give here an illustrative example of an $r|s$-form.
\begin{ex}
Let $\a^F\in V^*$ be an array of covectors of suitable parity. Then
from the properties of the Berezinian it follows that the
function
$L(v):=\Ber(\langle v_F, \a^G)\rangle)$ satisfies~\eqref{bers},\eqref{eqs}.
So it is a form. If $s>0$, $L$ has a pole at those odd vectors
whose linear span is not transverse to
the annihilator of the linear span of the odd part of $(\a^G)$. If $s=0$,
then $L(v)=\det(\langle v_i, \a^j)\rangle)$, where
$i,j=1,\dots,r$, so $L$ is nothing else than the exterior product
$\a^1\wedge\dots\wedge\a^r$. In general, this form with a
singularity should be regarded as a ``nonlinear analog" of the
exterior product of an array of even and odd covectors $\a^F$. It
naturally appears in physical context
(e.g.,~\cite{hov:bv},\cite{bel:pco}).
\end{ex}

As shown in~\cite{tv:dual}, the above construction of forms is not
sufficient and must be supplemented in order to obtain $r|s$-forms
with $r\in \mathbb Z$ arbitrary, including negative values. This
is achieved by the following ``dualization" and the subsequent ``stability
argument". When we shall need to distinguish forms in the sense of
Definition~\ref{straight}, we shall call them ``straight forms".
We shall denote the space of (straight) $r|s$-forms on $V$
by $\E^{r|s}(V)$.

\begin{de} A {\it dual form} on $V$ of codegree $p|q$ is a smooth
map $\L:\underbrace{V^*\times\dots V^*}_{p}\times
\underbrace{V^*\P\times\dots \times V^*\P}_{q} \to \Vol V$ satisfying
the
conditions
\begin{align}
  &\L(ph)=\L(p) \Ber h,\label{berr}\\
\intertext{for all $h\in \GL(p|q)$ and}
  &\dder{\L}{\pak}{\pbl}
  +(-1)^{\at\bt+(\at+\bt)\lt}\dder{\L}{\pbk}{\pal}=0\label{eq}.
\end{align} The arguments of $\L$ (covectors) are written
as vector-columns, and they are organized in a matrix $p=(\pak)$.
Notice that due to the condition~\eqref{berr}, odd covectors $p^K$, $\tilde
K=1$, should be linearly independent, hence $0\leq q\leq m$.
\end{de}

Fix a dimension $r|s$ and consider $V\oplus \R{r|s}$.
\begin{de}  A {\it mixed form} on $V$ of codegree $p|q$ and
additional degree $r|s$
is a smooth
map
$$\L:\underbrace{(V\oplus \R{r|s})^*\times\dots\times
(V\oplus \R{r|s})^*}_{p}\times
\underbrace{(V\oplus \R{r|s})^*\P\times\dots\times
(V\oplus \R{r|s})^*\P}_{q} \to \Vol V
$$ satisfying the following
conditions~(\ref{mberr})--(\ref{eq3}):
\begin{align}
  &\L(ph,wh)=\L(p,w) \Ber h,\label{mberr}\\
\intertext{for all $h\in \GL(p|q)$,}
&\L(p+aw,gw)=\L(p,w) \Ber g,\label{mberl}\\
\intertext{for all $g\in \GL(r|s)$ and all
$a\in\Mat(r|s\times n|m)$, and}
  &\dder{\L}{\pak}{\pbl}+
  (-1)^{\at\bt+(\at+\bt)\lt}\dder{\L}{\pbk}{\pal}=0,
\label{eq1}\\
&\dder{\L}{\pak}{\wfl}+
(-1)^{\at\ft+(\at+\ft)\lt}\dder{\L}{\wfk}{\pal}=0,
\label{eq2}\\
&\dder{\L}{\wfk}{\wgl}+
(-1)^{\ft\gt+(\ft+\gt)\lt}\dder{\L}{\wgk}{\wfl}=0
\label{eq3},
\end{align}
where $p=(\pak)$, $w=(\wfl)$ and for a given $K$ the entries
$\pak,\wfk$ are the components of a covector on $V\oplus
\R{r|s}$ (where $K$ is the number of the covector). Matrix notation
suggests placing $p$ over $w$ in the argument of $\L$, but for
typographic reasons we shall do it only when convenient. Notice
that $s\leq q \leq m+s$ because of \eqref{mberr},\eqref{mberl}.
\end{de}

Examples of dual and mixed forms can be mimicked from the examples of
straight forms (since they are defined via similar conditions),
and  we skip them.

Notation: $\E_{p|q}(V)$ and $\E_{p|q}^{r|s}(V)$ for the spaces of
dual and mixed forms on $V$, respectively. We shall omit the
indication to $V$ when no confusion is possible. Notice that
$\E_{p|q}(V)=\E_{p|q}^{0|0}(V)$

Consider the following homomorphisms:
$\s=\s_{k|l}: \E_{p|q}^{r|s}\to \E_{p+k|q+l}^{r+k|s+l}$ and
$\s^{-1}=\s_{k|l}^{-1}:\E_{p+k|q+l}^{r+k|s+l} \to \E_{p|q}^{r|s}$,

\begin{align}\label{stab}
  &(\s\L)\left(\begin{array}{cc}
         \pone & \ptwo  \\
         \woo & \wot\\
         \wto & \wtt
\end{array}\right):=\L\left(\begin{array}{c}
\pone-\ptwo\wtt^{-1}\wto \\
\woo-\wot\wtt^{-1}\wto \end{array}\right)\cdot\Ber\wtt,\\
  &(\s^{-1}\L^*)\left(\begin{array}{c}p \\w\end{array}\right):=\L^*
    \left(\begin{array}{cc}
         p & 0  \\
         w & 0\\
         0 & 1
        \end{array}\right),
\end{align}
where $\L\in \E_{p|q}^{r|s}$, $\L^*\in \E_{p+k|q+l}^{r+k|s+l}$.
(We write arguments of forms as matrices and subdivide them into
blocks corresponding to the ``first" and ``last" rows and columns.)

\begin{thm}[\cite{tv:dual}]
\label{stabil} Maps $\s$ and $\s^{-1}$ are well-defined (in
particular, $\s$ uniquely extends to all admissible arguments of
$\L$) and are indeed mutually inverse isomorphisms of the spaces
$\E_{p|q}^{r|s}$ and $\E_{p+k|q+l}^{r+k|s+l}$. The equality
$\s_{k|l}\s_{k'|l'}=\s_{k+k'|l+l'}$ holds.
\end{thm}

Define $\EE^{k|l}(V):=\dlim \E_{p+N|q+M}^{r+N|s+M}(V)$, where
$k|l=r+n-p|s+m-q$ and call it the space of {\it stable
$k|l$-forms} on $V$. Note that $k\in \mathbb Z$ (may be negative),
while $l=0,\dots,m$ . It's not hard to produce an example of a
stable $k|l$-form with negative $k$ (if $l > 0$). Take as a
representative a dual form with the number of even arguments
greater that $n$ (exactly as in  examples of straight $r|s$-forms
with $r>n$, cf.~\cite{tv:git}). Similarly, if $l < m$, there are
nonzero $k|l$-forms with $k>n$.

Obviously, $\EE^{k|l}(V)\cong\E_{p|q}^{r|s}(V)$ if $k=r+n-p$,
$l=s+m-q$, for all $r,s,p\geq 0$ and $s\leq q\leq
s+m$.

\begin{co}$  \EE^{k|l}(V)\cong
\E_{n-k|m-l}(V)$ for $k\leq n$.
\end{co}

Consider the following homomorphisms: $\t=\t_{r|s}:\E^{r|s}\to
\E_{n|m}^{r|s}$ and $\t^{-1}=\t_{r|s}^{-1}:\E_{n|m}^{r|s}\to \E^{r|s}$,
\begin{align}\label{iso}
  (\t L)
  \left(
  \begin{array}{c}p \\w
  \end{array}
  \right)&:=L(wp^{-1})\cdot\Ber p,\\
  (\t^{-1} \L)(v)&:=\L
  \left(
    \begin{array}{c}1 \\
     v
    \end{array}
  \right),
\end{align}
where $\L\in \E_{n|m}^{r|s}$, $L\in \E^{r|s}$.

\begin{thm}[\cite{tv:dual}]
\label{isom}Maps $\t$ and $\t^{-1}$ are well-defined (in
particular, $\t$ uniquely extends to all admissible arguments of
$L$) and are indeed mutually inverse isomorphisms of the spaces
$\E_{n|m}^{r|s}$ and $\E^{r|s}$.
\end{thm}

\begin{co}$\EE^{k|l}(V)\cong \E^{k|l}(V)$  for
$k\geq 0$.
\end{co}

\begin{rem}In view of  Theorems~\ref{stabil} and \ref{isom} one may
regard it excessive to consider all spaces of mixed forms. Indeed, it is
sufficient to consider only $\E^{r|s}$ and $\E_{p|q}$ together
with the isomorphism $\E^{r|s}\cong\E_{n-r|m-s}$ defined in the range
$0\leq r\leq n$. However, it would be
practically restrictive. It is easier to work with various
operations in terms of mixed forms.
\end{rem}

\subsection{Operators $e(\a)$, $e(u)$.
Commutation relations and the module structure}

Consider a covector $\a\in V^*$. We introduce an operator
$e(\a):\E^{r|s}_{p|q}\to \E^{r+1|s}_{p|q}$  by the following formula:
\begin{equation}\label{ext}
  e(\a)\L:= (-1)^r \a_A w_{r+1}^K (-1)^{\tilde \a\tilde
  A}\der{\L}{\pak},
\end{equation}
where $\a=e^A\a_A$.

Likewise, consider a vector $u\in V$. Define
$e(u):\E^{r|s}_{p|q}\to \E^{r|s}_{p+1|q}$ by the formula
\begin{multline}\label{int}
  e(u)\L:= \\
  (-1)^r u^A \left(
  p_A^{p+1} - (-1)^{\bt\kt}
  \pak p_B^{p+1}\der{}{\pbk}-(-1)^{\ft\kt}\pak w_F^{p+1}\der{}{\wfk}
  \right)\L,
\end{multline}
where $u=u^A e_A$. Here $(e_A)$ and $(e^A)$ are dual bases of $V$ and
$V^*$.

\begin{rem}
On dual forms, $e(u):\E_{p|q}\to\E_{p+1|q}$,
\begin{equation}\label{intdual}
  e(u)\L=
  (-1)^r u^A \left(p_A^{p+1} - (-1)^{\bt\kt}  \pak
  p_B^{p+1}\der{}{\pbk}\right)\L.
\end{equation}

\end{rem}

The proof that $e(\a)$ and $e(u)$ indeed map forms to forms and
do not depend on the choice of bases  is postponed until
Section~\ref{analiz}.
The
parities of $e(\a)$ and $e(u)$ are the same as the respective
parities of $\a$ and $u$; operators $e(\a)$ and $e(u)$ depend on $\a$ and
$u$ linearly.

\begin{thm}
\label{stability}
The operators $e(\a)$ and $e(u)$ are stable, i.e., they
commute with all isomorphisms $\s_{k|l}$.  Under the
isomorphism~\eqref{iso}, the operator $e(\a)$ corresponds to the operator
$e_{\a}: \E^{r|s}\to \E^{r+1|s}$,
\begin{align}\label{ext1}
    e_{\a}&=(-1)^r \left(
    \vaone\a_A -(-1)^{\aat\ft+\bt}\vfa\a_A\,\vbone\,\der{}{\vfb}
           \right)
\intertext{and if
$r>0$ the operator $e(u)$ corresponds to
the operator $i_u: \E^{r|s}\to \E^{r-1|s}$,}
\label{int1}
  i_u&=(-1)^{r-1}u^A\,\der{}{\vra},
\end{align}
the substitution of $u\in V$ into the last even slot of
$L\in\E^{r|s}$. Here
$L=L(v)$, $v=(\vfa)$. (The operators $e_{\a}$, $i_u$ were introduced
in~{\em \cite{tv:git}}.)
\end{thm}
\begin{proof} Consider $e(u)$. We have
to check that $e(u)$ commutes with $\s_{1|0}$ and $\s_{0|1}$.
We shall consider $\s_{1|0}$ (the case of $\s_{0|1}$ is similar, but
simpler). Denote $\s:=\s_{1|0}$. It is sufficient to give proof
for $\L\in\E_{p|q}$, then the general case will follow. Consider the
diagram
\begin{equation}
\begin{CD}
    \E_{p|q}@>{\s}>>\E_{p+1|q}^{1|0}\\
    @V{e(u)}VV        @VV{e(u)}V\\
    \E_{p+1|q}@>>{\s}>\E_{p+2|q}^{1|0}
\end{CD}
\end{equation}
Take $\L\in\E_{p|q}$. Apply $\s$. We get $\L^*\in\E_{p+1|q}^{1|0}$,
where $\L^*\left(\begin{smallmatrix}p & p'\\
w & w'
\end{smallmatrix}\right)=\L\left(p-p'{w'}^{-1}w\right)\,w'$. Here $p=(\pak)$,
$w=(\wfk)$,
$p'=(\paone)$, $w'=w^{p+1}$. Apply $e(u)$. We obtain
\begin{multline}\label{sigmae}
    (e(u)\L^*)\left(
        \begin{matrix}p & p' &p''\\
                           w & w'&w''
        \end{matrix}\right)=
        -u^A\left(\patwo
        -(-1)^{\bt{\tilde{K^*}}}\pakstar\pbtwo\,\der{}{\pbkstar}\,- \right.
        \\
        \left.
        \pakstar\wtwo\,\der{}{\wkstar}
        \right)\L^*=
        -u^A\left(
        \patwo -(-1)^{\bt\kt}\pak\pbtwo\,\der{}{\pbk}-
        \right.
        \\
        \left.
        \paone\pbtwo\,\der{}{\pbone}\,-
        \pak\wtwo\,\der{}{\wk} -
         \paone\wtwo\,\der{}{\wone}
        \right)
         \L\left(p-p'{w'}^{-1}w\right)\,\wone
         \\=
        -u^A\left(\patwo\L \wone -
        (-1)^{\bt\kt}\pak\pbtwo\,\der{\L}{\pbk} \wone +
        \right.\\
        \paone\pbtwo\wk\,\der{\L}{\pbk}  +
        \pak\wtwo\pbone\,
         \der{\L}{\pbk} (-1)^{\bt\kt}+
        \paone\wtwo\pbone\wk\,\der{\L}{\pbk}
        \\
        \left.\bigl(-\frac{1}{(\wone)^2}\bigr)\wone
        - \paone\wtwo\,\L \right),
\end{multline}
where in the last expression the argument of $\L$ and $\lder{\L}{p}$ is
$p-p'{w'}^{-1}w$ and we denote $p'':=(\patwo)$, $w'':=(\wtwo)$.
Now let us apply first $e(u)$, then $\s$. Calculate:
\begin{equation}
  (e(u)\L)\bigl(\begin{matrix}p &p''\end{matrix}\bigr)=
  u^A\left(\patwo-(-1)^{\bt\kt}\pak\pbtwo\,\der{}{\pbk}\right)\L(p);
\end{equation}
applying $\s$ we obtain
\begin{multline}
(\s e(u)\L)\left(
        \begin{matrix}p & p'' &p'\\
                      w & w''&w'
        \end{matrix}\right)=\left(e(u)\L\right)
        \bigl(p-p'{w'}^{-1}w,
        p''-p'(\wone)^{-1}\wtwo\bigr)\,\wtwo
        \\=u^A\Biggl( (\patwo - \paone(\wone)^{-1}\wtwo) \, \L
        - (-1)^{\bt\kt} (\pak - \paone(\wone)^{-1}\wk) \Biggr.
        \\ \Biggl.
        (\pbtwo - \pbone(\wone)^{-1}\wtwo)
        \der{\L}{\pbk}\Biggr) \wone,
\end{multline}
where the argument of $\L$ and $\lder{\L}{p}$ in the last expression
is $p-p'{w'}^{-1}w$.
Multiplying through, we obtain exactly the same terms as
in~\eqref{sigmae}  with the opposite sign. Notice that $\s e(u)\L$ as a
form is skew-symmetric in even columns. Thus we can swap
$\left(\begin{smallmatrix}p' \\ w'\end{smallmatrix}\right)$ and
$\left(\begin{smallmatrix}p'' \\ w''\end{smallmatrix}\right)$, cancelling
the minus sign,
and obtain
\begin{equation}
  (\s e(u)\L)\left(
        \begin{matrix}p & p' &p''\\
                      w & w'&w''
        \end{matrix}\right)=
  (e(u)\s \L)\left(
        \begin{matrix}p & p' &p''\\
                      w & w'&w''
        \end{matrix}\right),
\end{equation}
as desired. Stability of $e(\a)$ is proved in the same way, and we
omit the calculation.

Let us turn to the relation with the isomorphisms~\eqref{iso}.
Consider  the following diagram.

\begin{diagram}
\E_{n|m}^{r|s} & \rTo^{e(u)} & \E_{n+1|m}^{r|s} & \rTo^{\s^{-1}} & \E_{n|m}^{r-1|s}\\
\dTo<{\t^{-1}}           &      &                  &       &   \dTo >{\t^{-1}} \\
\E^{r|s}       &         &     \rTo_{i_u}      &       &\E^{r-1|s}   \\
\end{diagram}
The claim is that it is commutative. To check this, take $\L\in
\E_{n|m}^{r|s}$. We have:
\begin{multline*}
  (i_u\t^{-1}\L)(v)=(-1)^{r-1}\ua \,\der{}{\vra}(\t^{-1}\L)(v)
  =(-1)^{r-1}\ua\,\der{}{\vra}\L
  \begin{pmatrix}
      1 \\
      v
  \end{pmatrix}=
  \\ (-1)^{r-1}\ua \,\der{\L}{\wra}
  \begin{pmatrix}
      1 \\
      v
  \end{pmatrix};
\end{multline*}
now,
\begin{align*}
    \begin{split}
    (e(u)\L)
     \begin{pmatrix}
     p & \pnone \\
     w  & \wnone
     \end{pmatrix}=(-1)^r \ua\left(\panone-(-1)^{\bt\kt}\pak\pbnone\,
    \der{}{\pbk} - \right.
     \\
    \left.(-1)^{\ft\kt}\pak\wfnone\,\der{}{\wfk}\right)
     \L\begin{pmatrix}
      p \\
      w
     \end{pmatrix};
     \end{split}\\
    \begin{split}
     &(\s^{-1}e(u)\L)\begin{pmatrix}p \\w^*\end{pmatrix}=
     (e(u)\L)\begin{pmatrix}p & \pnone \\w  & \wnone
             \end{pmatrix}_{\left|\begin{aligned}&\scriptstyle
                                 w_r^{n+1}=1\\
                                 &\scriptstyle
                                 \wrk  =0 \quad (K\neq n+1)\\
                                 &\scriptstyle
                                 \wfnone =0 \quad (F\neq r)\\
                                 &\scriptstyle
                                 \panone=0
                       \end{aligned}
               \right.}=
  \\&(-1)^{r}u^A\left(0-(-1)^{0}\pak\,\der{}{\wrk}\right)
  \L\begin{pmatrix}p\\
  w^*\\
  0 \end{pmatrix}=
  (-1)^{r}u^A\left(-\pak\,\der{\L}{\wrk}\begin{pmatrix}p\\
                                                w
                                          \end{pmatrix}
             \right);
\end{split}\\
\intertext{hence}
&(\t^{-1}\s^{-1}e(u)\L)(v)=(-1)^r \left(-u^A\der{\L}{\wra}
  \begin{pmatrix}
    1 \\
    v
  \end{pmatrix}
\right)=i_u \t^{-1}\L(v),
\end{align*}
as desired. (Here  $w^*$ stands for $w$ without the row $w_r$.)
In a similar way the equality $e(\a)\t=\t e_{\a}:\E^{r|s}\to\E_{n|m}^{r+1|s}$
is checked.
\end{proof}

\begin{co}For exterior forms on a purely even space $V$ the operator
$e(\a)$ corresponds to the usual
exterior
multiplication $\a\,\wedge\,$.
The operator $e(u)$ corresponds to the usual interior multiplication or
contraction  $i_u=u\cowedge\,$.
\end{co}

Note that in our mixed description both operators increase respective
degrees and thus have appearance of ``exterior" products.

\begin{thm} \label{commut}The operators $e(\a)$ and $e(u)$ obey the following
relations:
\begin{align}
  e(u) e(v) + (-1)^{\ut\vt} e(v) e(u)  & =0, \label{commut1}\\
  e(\a) e(\b) + (-1)^{\tilde\a\tilde\b}e(\b) e(\a)  & =0, \label{commut2}\\
  e(u) e(\a) + (-1)^{\tilde\a\ut}e(\a) e(u)  & = \langle
  u,\a\rangle\,\s.\label{cliff}
\end{align}
Here $u,v\in V$, $\a,\b\in V^*$, and
$\s=\s_{1|0}:\E_{p|q}^{r|s}\to \E_{p+1|q}^{r+1|s}$ is the
stability isomorphism~\eqref{stab}.
\end{thm}
\begin{proof}
To find relations between $e(u)$ and $e(v)$, for $u,v\in V$, it is
sufficient to consider the case $r=s=0$. (The general case is
formally reduced to it by considering dual forms on extended space
$V\oplus \R{r|s}$ and by setting $u^F=v^F=0$.) Then for
$\L\in\E_{p|q}$ we have
\begin{multline}\label{vykl1}
e(u)\,e(v)\L
\\=u^A\left(p_A^{p+2}-(-1)^{\bt\kt}\pak
p_B^{p+2}\,\der{}{\pbk} \right)
v^C\left(p_C^{p+1}-(-1)^{\dt\lt}\pcl
p_D^{p+1}\,\der{}{\pdl} \right)\L
\\=u^A v^C (-1)^{(\vt+\ct)\at}
\Biggl(\patwo\pcone - \paone\pctwo -
(-1)^{\ct\dt}\patwo\pdone\pcl\,\der{}{\pdl}-\\
(-1)^{\bt\ct+\at(\bt+\ct)}\pbtwo\pcone\pal\,\der{}{\pbl} +
(-1)^{\at(\ct+\dt)}\pctwo\pdone\pal\,\der{}{\pdl}+\\
(-1)^{\ct\dt}\paone\pdtwo\pcl\,\der{}{\pdl}+
 (-1)^{a}\,
\pbtwo\pdone  \pak\pcl\,\dder{}{\pbk}{\pdl}\Biggr)\L,
\end{multline}
where $a=\bt\ct+\bt\lt+\bt\dt+\ct\kt+\kt\lt+\at\bt+\at\dt+\ct\dt$.
Notice that the range of $K$ in
the first line of~\eqref{vykl1} contains $p+1$.
Simultaneously interchanging $u$ and $v$ and  the indices $A$ and
$C$,
we obtain
\begin{multline}
e(v)\,e(u)\L
\\= (-1)^{\ut\vt}u^A v^C (-1)^{(\vt+\ct)\at}
\Biggl(\paone\pctwo - \patwo\pcone -
(-1)^{\at(\ct+\dt)}\pctwo\pdone\pal\,\der{}{\pdl}
\\ -
(-1)^{\ct\dt}\paone\pdtwo\pcl\,\der{}{\pdl}+
(-1)^{\ct\dt}\patwo\pdone\pcl\,\der{}{\pdl} +
\\(-1)^{\at\dt+\at\ct+\ct\dt}\pdtwo\pcone\pal\,\der{}{\pdl}+
(-1)^b\pbtwo\pdone\pak\pcl\,\dder{}{\pdk}{\pbl}
\Biggr)\L,
\end{multline}
where $b=\ct\kt+\at\bt+\kt\lt+\bt\ct+\ct\dt+\at\dt+\lt\dt$. Now we
see that all terms  except for the
last one in $(-1)^{\ut\vt}e(v)e(u)\L$
would cancel the similar terms in $e(u)e(v)\L$.
Notice that $a+b=\bt\dt+(\bt+\dt)\lt$. It follows that
\begin{multline}
    \left(e(u)e(v)+(-1)^{\ut\vt}e(v)e(u)\right)\L=(-1)^a \pbtwo\pdone\pak\pcl\,
   \\ \left( \dder{\L}{\pbk}{\pdl}+
   (-1)^{\bt\dt+(\bt+\dt)\lt}\dder{}{\pdk}{\pbl}\right),
\end{multline}
which equals zero by the equation~\eqref{eq}.

Consider now $e(\a)$ and $e(\b)$. For $\L\in\E_{p|q}^{r|s}$ we
readily
have
\begin{multline}\label{vykl2}
    e(\a)e(\b)\L = (-1)^{r+1}\a_A \wktwo \,\der{}{\pak}\,
    \left(
    (-1)^r\b_B
    \wltwo\,\der{\L}{\pbl}(-1)^{\tilde\b \bt}
    \right)=
    \\-(-1)^{\aat\at+\bbt\bt}\a_A\b_B\wktwo\wlone\,\dder{\L}{\pak}{\pbl}
    (-1)^{(\bbt+\bt)\at+(\at+\kt)\lt}.
\end{multline}
Similarly, for $e(\b)e(\a)$ we obtain
\begin{multline}
    e(\b)e(\a)\L = -(-1)^{\aat\bbt +\aat\at +\bbt\bt +(\bt+\kt)\lt +\at\bbt}
    \a_A\b_B\wktwo\wlone\,\dder{\L}{\pbk}{\pal}=
    \\
    (-1)^{\aat\at+\bbt\at_\aat\bbt+\kt\lt+\at\bt+\at\lt}\a_A\b_B\wktwo\wlone\,
    \dder{\L}{\pak}{\pbl}=
    \\-(-1)^{\aat\bbt}e(\a)e(\b)\L,
\end{multline}
again by the equation~\eqref{eq}.

Finally, let us find the relation between operators $e(u)$ and
$e(\a)$. Notice that
$e(u)e(\a), \ e(\a)e(u):\E_{p|q}^{r|s}\to\E_{p+1|q}^{r+1|s}$. For
$\L\in\E_{p|q}^{r|s}$ by a direct calculation similar to~\eqref{vykl1},\eqref{vykl2}
using the equations~\eqref{eq1},\eqref{eq2}, we obtain the
equality
\begin{multline}\label{cliff1}
    \left(e(u) e(\a) + (-1)^{\tilde\a\ut}e(\a) e(u)\right)\L   =
    \\u^A \a_A \left(
    w_{r+1}^{p+1}
    -(-1)^{\bt\kt}\wkone\pbone\,\der{}{\pbk}-(-1)^{\ft\kt}\wkone
    w_{F}^{p+1}\,\der{}{\wfk}
    \right)\L.
\end{multline}
Apply now the transformation $\s^{-1}:\E_{p+1|q}^{r+1|s}\to
\E_{p|q}^{r|s}$. That means setting $w_{r+1}^{p+1}:=1$,
$\wkone:=0$, $\pbone:=0$, $w_{F}^{p+1}:=0$. We arrive at
\begin{equation}
  \s^{-1}\left(e(u) e(\a) + (-1)^{\tilde\a\ut}e(\a) e(u)\right)\L=\langle
  u,\a\rangle\L,
\end{equation}
from where~\eqref{cliff} follows. Notice that by this calculation
we showed that the operator in the r.h.s. of~\eqref{cliff1} gives
another expression for the isomorphism
$\s_{1|0}:\E_{p|q}^{r|s}\to\E_{p+1|q}^{r+1|s}$.
\end{proof}

\begin{co} {\em(1)} The space $\E_{\bdot|q}^{\bdot|s}(V)$ is a module over
exterior algebras $\E^{\bdot}(V)$ and $\E^{\bdot}(V^*)$ defined by
relations $u v=-(-1)^{\ut\vt}v u$ and $\a \b=-(-1)^{\tilde
\a\tilde \b}\b \a$.

{\em(2)} The space of stable forms $\EE^{\bdot|s}(V)$ is a module over
a Clifford algebra $\Cliff(V\oplus V^*)$  defined by relations $u
v=- (-1)^{\ut\vt}v u$, $\a \b=- (-1)^{\tilde \a\tilde \b}\b \a$
and $u\a +(-1)^{\ut\tilde \a}\a u=\langle u,\a\rangle$.
\end{co}

\begin{rem}
Notice that we arrive at the relations
of exterior and Clifford algebras (in ``skew" versions) not as
conventions but as  actual identities between linear operators.
It is also worth  noting that the anticommutation relations
obtained here for $e(u)$ and $e(\a)$ are not at all obvious.
While under the isomorphism with straight or dual
forms one of the operators
$e(u)$ or $e(\a)$  can be interpreted as a substitution into a suitable even
slot (hence the anticommutativity between such operators will become transparent),
the other one will remain an  ``exterior product"
defined by a  formula like~\eqref{ext1}, which involves both even and odd
slots. By duality $e(u)$ transforms into $e(\a)$ and vice versa. However,
this can be
exploited only in the common range $0\leq r\leq n$ where dual and
straight forms are both good. Hence a  certain portion of tedious
calculations is unavoidable to get all the
relations~{\eqref{commut1}--\eqref{cliff}}.
\end{rem}

\section{Cartan calculus}\label{analiz}

\subsection{Differential}
Consider a supermanifold $M=M^{n|m}$. For forms on $M$, i.e.,
sections of the corresponding vector bundles associated with $TM$,
we shall use the notation $\O^{r|s}$, $\O_{p|q}$, $\O^{r|s}_{p|q}$
and $\boldsymbol \O^{r|s}$. By $\O^{\bdot}=\oplus\O^{k}$ we shall
denote the algebra of ``naive" differential forms with the
skew-commutative convention (and the even differential,
cf.~\cite{manin}). A differential
$\D:\O^{r|s}_{p|q}\to\O^{r+1|s}_{p|q}$ is defined by the formula
\begin{equation}\label{dif}
  \D\L:=(-1)^r \wkone(-1)^{\at\kt}\,\der{}{\xa}\der{\L}{\pak}
\end{equation}
(see~\cite{tv:dual}). In~\cite{tv:dual} it is proved that the operator
$\D$ is stable, hence we have a complex
$\D:\OO^{\bdot\,|s}\to\OO^{\bdot\,+1|s}$. For $\bdot\geq 0$, this
complex is isomorphic to the ``straight" complex
$d:\O^{\bdot\,|s}\to\O^{\bdot\,+1|s}$ studied
in~\cite{tv:git} and for $\bdot\!\leq n$ to
the complex of dual forms
$\DD:\O_{n-\bdot+1|m-s}\to\O_{n-\bdot|m-s}$ introduced
in~\cite{tv:dual}:

{
\newcommand{\snul}{\O^{0|s}}
\newcommand{\sodin}{\O^{1|s}}
\newcommand{\sn}{\O^{n|s}}
\newcommand{\snodin}{\O^{n+1|s}}
\newcommand{\omin}{\OO^{-1|s}}
\newcommand{\onul}{\OO^{0|s}}
\newcommand{\oodin}{\OO^{1|s}}
\newcommand{\on}{\OO^{n|s}}
\newcommand{\onodin}{\OO^{n+1|s}}
\newcommand{\dnodin}{\O_{n+1|m-s}}
\newcommand{\dn}{\O_{n|m-s}}
\newcommand{\dmin}{\O_{n-1|m-s}}
\newcommand{\dnul}{\O_{0|m-s}}

\newarrow{Isom}=====

\begin{diagram}[width=2em,height=1.5em]
 & & 0&\rTo
 &\snul&\rTo&\sodin&\rTo&\dots&\rTo&\sn&\rTo&\snodin&\rTo&\dots\\
 & & & & \dIsom & & \dIsom & & & & \dIsom & & \dIsom& & \\
 \dots&\rTo&\omin&\rTo&\onul&\rTo&\oodin&\rTo&\dots&\rTo&\on&\rTo&\onodin&\rTo&\dots\\
  & &\dIsom & &\dIsom& &\dIsom& & & & \dIsom & & & & \\
  \dots&\rTo&\dnodin&\rTo&\dn&\rTo&\dmin&\rTo&\dots&\rTo&\dnul&
  \rTo&0
  & & \\
\end{diagram}
}
(vertical lines are isomorphisms).

Consider a mixed form $\L$ and a function $f$. Calculate
$\D(f\L)$:
\begin{multline}
    \D(f\L)=(-1)^r \wkone(-1)^{\at\kt}\,\der{}{\xa}\der{}{\pak}(f\L)
    =\\(-1)^r \wkone(-1)^{\at\kt}\,\der{}{\xa}f\der{\L}{\pak}(-1)^{\ft(\at+\kt)}
    =\\(-1)^r \wkone(-1)^{\at\kt}\left((-1)^{\ft(\at+\kt)}\da f\,\der{\L}{\pak}
    +(-1)^{\ft\kt}f\,\der{}{\xa}\der{\L}{\pak}\right)
    =\\ f\,\D\L + (-1)^r\da\,f\wkone\,\der{\L}{\pak}(-1)^{\ft\at}
    =f\,\D\L + e(df)\,\L,
\end{multline}
where $df=dx^A\da f$ is considered as an element of $\O^1(M)$. We
stress that the algebra with the {\it even} differential is
considered. Since $\D(f\L)$ is a form and $f\,\D\L$ is a form,
it follows that $e(df)\,\L$ is a well-defined form.
We can conclude that
for arbitrary $1$-form $\a$ the operation $e(\a)$ is also well-defined, i.e.,
does not
depend on the choice of coordinates and maps mixed forms into
mixed forms.
The formula~\eqref{ext} is extracted from this calculation.
Similar calculation  gives the formula~\eqref{ext1} for $e_{\a}$
on straight forms; by duality it can be rewritten to produce a
formula~\eqref{intdual} for $e(u)$ on dual forms,
from which we get our formula~\eqref{int}
on mixed forms. Thus it follows that both operators $e(u)$,
$e(\a)$ on mixed forms
are well-defined, which justifies our consideration in the previous section.
It is not easy to give a purely algebraic proof of this fact.

\begin{rem}
    The stability of $e(u)$, $e(\a)$ as well can be deduced from the
    stability of $\D$.
\end{rem}

In the previous Section we got the module structure of mixed forms
over $\O^{\bdot}(M)$.

\begin{thm}
\label{leibniz}Leibniz formula holds:
\begin{equation}\label{leib}
  \D(\o\,\L)= d\o\,\L + (-1)^k \o\,\D\L,
\end{equation}
for $\o\in\O^k$ and $\L\in\O^{r|s}_{p|q}$.
\end{thm}
\begin{proof}
    Since $\O^{\bdot}(M)$ is a differential graded algebra,
    generated by elements $df$ over $C^{\infty}(M)$ (locally),
    it is sufficient to check the formula~\eqref{leib} for two
    cases: $\o=f$ and $\o=df$, where $f$ is a function. The first
    case was considered above. Consider $\o=df$. Then, by
    definition,
    \begin{equation}
    df\,\L=\D(f\L) - f\,\D\L.
    \end{equation}
    Apply $\D$. We obtain
    \begin{multline}
    \D(df\,\L)=\D\D(f\L) - \D(f\,\D\L)
    =0-df\,\D\L - f\D\D\L=-df\,\D\L
    =\\ ddf\,\L+(-1)^1 df\,\D\L,
    \end{multline}
    as desired.
\end{proof}
Therefore, $\OO^{\bdot\,|s}$ is a graded differential module over
$\O^{\bdot}$ for all s.

\begin{rem} Notice that $\E^{\bdot}\cong\EE^{\bdot|0}$,
$\O^{\bdot}\cong\OO^{\bdot|0}$   as
modules.
\end{rem}

\subsection{Homotopy identity}
Consider a vector field $X\in \Vect M$ and the corresponding
infinitesimal transformation: $x^A\mapsto x^A + \e X^A(x)$,
$\e^2=0$. By a straightforward calculation we obtain the following
formula for the Lie derivative on mixed forms:
\begin{equation}
    \d_X\L=
    X^A\,\der{\L}{x^A}- (-1)^{\at\xt}\der{X^B}{x^A}\
    \pbk\der{\L}{\pak}     + (-1)^{\at(\xt+1)} \der{X^A}{x^A}\,\L,
\end{equation}
where we picked the notation $\d_X$ to avoid overloading the
letter {\it `L'}. The Lie derivative $\d_X$ has the same parity
as $X$. It preserves all degrees and is obviously a derivation for
all kinds of natural multiplications. Operation $\d_X$ commutes
with the stability isomorphisms~\eqref{stab} and with the
isomorphisms~\eqref{iso}.

\begin{thm} \label{cart}
For mixed forms on a supermanifold $M$, the following identity holds:
\begin{equation}\label{cartan1}
  \D\, e(X) + e(X)\,\D = \d_X\,\s,
\end{equation}
where $\s=\s_{1|0}:\O^{r|s}_{p|q}\to\O^{r+1|s}_{p+1|q}$ is the
stability isomorphism.
\end{thm}
\begin{proof}
    Let $\L$ be in $\O_{p|q}^{r|s}$. Consider
    $\s^{-1}:\O_{p+1|q}^{r+1|s}\to\O_{p|q}^{r|s}$. Recall that the
    action of this operator consists in setting $\paone=0$, $w_{F}^{p+1}=0$,
    $\wkone=0$, $w_{r+1}^{p+1}=1$ in the argument. We shall find
    $\s^{-1}e(X)\D\L$ and $\s^{-1}\D e(X)\L$. Directly from~\eqref{int}:
    \begin{multline}\label{ed}
        \s^{-1}e(X)\D\L=(-1)^{r+1}\Xa\left(-\pak\,\der{}{\wkone}\,\D\L\right)
        =\\
        (-1)^{r}\Xa \pak\,\der{}{\wkone}\,
        \left((-1)^r\wlone(-1)^{\bt\lt}\der{}{\xb}\der{}{\pbl}\L\right)
        =\\
        \Xa\pak(-1)^{\bt\kt}\der{}{\xb}\der{\L}{\pbk};
    \end{multline}
    now,
    \begin{multline}
        \s^{-1}\D e(X)\L=(-1)^r
        w_{r+1}^{K^*}(-1)^{\at\tilde{K^*}}\der{}{\xa}\der{}{\pakstar}
        (e(X)\L)_{\left|
               \begin{aligned}
               \scriptstyle w_F^{p+1}&\scriptstyle=0,\quad
               & \scriptstyle\paone&\scriptstyle=0, \\
               \scriptstyle\wkone&\scriptstyle=0,\quad
               & \scriptstyle w_{r+1}^{p+1}&\scriptstyle=1
               \end{aligned}
        \right.}
        =\\
        (-1)^r\left(\der{}{\xa}\der{}{\paone}(e(X)\L)\right)_{\left|
        \scriptstyle p^{p+1}=0, \ w^{p+1}=0\right.}
        =\left(\der{}{\xb}\der{}{\pbone}\left(
        \Xa\Bigl(p_A^{p+1}\L
        \Bigr.\right.\right.
        \\
        \Bigl.\left.\left.
        -(-1)^{\ct\kt}\pak p_C^{p+1}\,\der{\L}{\pck}
        -(-1)^{\ft\kt}\pak w_F^{p+1}\,\der{\L}{\wfk}
           \Bigr)
                                        \right)
        \right)_{\left|
        \scriptstyle p^{p+1}=0,\  w^{p+1}=0\right.}
        =\\
        \der{}{\xb}\left(
            \Xa(-1)^{\bt(\at+\xt)}\left(\dab\L-(-1)^{\at\bt}
            \pak\,\der{\L}{\pbk}
            \right)
        \right)
        =\\
        (-1)^{\bt(\xt+1)}\der{\Xb}{\xb}\,\L + \Xb\,\der{\L}{\xb}\,
        -\der{\Xa}{\xb}(-1)^{\bt\xt}\pak\,\der{\L}{\pbk}\,
        \\
        -(-1)^{\bt\kt}\Xa \pak \,\der{}{\xb}\der{\L}{\pbk}.
    \end{multline}
    Comparing with~\eqref{ed}, we immediately conclude that
    \begin{multline}\label{sum}
        \s^{-1}\bigl(e(X)\,\D + \D \,e(X)\bigr)\L=
        \\
        (-1)^{\bt(\xt+1)}\der{\Xb}{\xb}\,\L + \Xb\,\der{\L}{\xb}\,
        -\der{\Xa}{\xb}(-1)^{\bt\xt}\pak\,\der{\L}{\pbk}=\d_X \L.
    \end{multline}
    Applying $\s$ to both sides of~\eqref{sum}, we obtain the
    desired
    identity~\eqref{cartan1}. (Notice that $\s$ and $\d_X$ commute.)
\end{proof}

\begin{co}\label{ccart}
In the complex of stable forms $\OO^{\bdot|s}$ we
have the usual form of  ``Cartan's homotopy identity":
\begin{equation}\label{cartan2}
  \D \,e(X) + e(X)\,\D = \d_X.
\end{equation}
\end{co}

\section{Discussion}\label{discuss}

We introduced the operators $e(u)$ and $e(\a)$ on the space of
mixed forms, where $u$ is a vector and $\a$ is a covector. They
are analogs of the contraction $u\cowedge\,$ and of the exterior product $\a\wedge\,$
for
usual forms on purely even vector space. Though these operations
change only even part of
degrees, their construction involves all (even and odd) arguments. We proved
that these operations are stable, hence they induce the
corresponding operations on the space of stable forms. We
established the anticommutation relations for the operators $e(u)$ and
$e(\a)$. They yield the relations of a super Clifford algebra (or, before stabilization,
with an additional central element
$\s$). It is remarkable that a ``skew-commutative" version of Clifford
relations (anticommutators without parity reversion) rather than more
popular choice of
commutators and reversed
parity naturally appears here.

The main incentive of considering these operators was the
necessity to straighten out the Cartan calculus for forms on
supermanifolds. The homotopy identity found in~\cite{tv:git} was
valid only for $r|s$-forms with $r>0$; the case $r=0$ had to be
mended with the help of an {\it ad hoc} augmentation. The
existence of Bernstein-Leites integral forms of negative degree
has given another hint to a ``hidden" part of the super
{Cartan-de~Rham}
complex. This hidden part was discovered in~\cite{tv:dual}. The
entire complex (incorporating positive and negative halves)
is made up by stable forms, for which mixed forms are
representatives. In the current paper we established the relation
between the differential and the operator $e(X)$, where $X$ is a
vector field. Again, for mixed forms it contains the element $\s$
and after stabilization an analog of the usual form of the homotopy identity is
reproduced. Thus, the introduction of the stable complex indeed
solves the problem.

What is next? We need to check the functorial behaviour of
stable forms and get a ``generalized" version of the homotopy
identity, which will imply the homotopy invariance of the complex
(note that $\d_X$ in~(\ref{cartan1},\ref{cartan2})
corresponds to an infinitesimal diffeomorphism;
we need perturbations of arbitrary maps),
hence
an analog of the Atiyah-Hirzebruch
sequence (cf.~\cite{tv:git}). The investigation of ``point cohomology" of stable forms
will
require more detailed analysis of their algebraic properties.
Another topic, which we did not touch here at all, is, of course,
integration. We hope to consider these subjects elsewhere.
In the paper~\cite{tv:lag}, the author showed that the variational
differential can be used to make a complex of arbitrary
Lagrangians of paths, not just forms. It would be interesting to
combine this fact with the results of~\cite{tv:dual} and of the
current paper.

\end{document}